\documentclass[10pt]{amsart}
\usepackage[utf8]{inputenc}
\usepackage[usenames, dvipsnames]{color}
\usepackage{mathtools}
\usepackage{amscd}
\usepackage{amsmath}
\usepackage{amssymb}
\usepackage{amsfonts}
\usepackage{amsthm}
\usepackage{enumitem}
\usepackage{stmaryrd}
\usepackage[notcite,notref]{showkeys}
\usepackage{soul}
\input xy

\newcommand{\ZZ}{{\mathbb{Z}}}

\newcommand\scalemath[2]{\scalebox{#1}{\mbox{\ensuremath{\displaystyle #2}}}}

\usepackage[all]{xy}
\usepackage{tikz-cd}
\usepackage{caption}
\usepackage{hyperref}

\newcommand{\CC}{{\mathcal C}}

\newcommand{\Z}{{\mathbb Z}}

\newcommand{\cB}{\mathcal{B}}

\usepackage{dsfont}

\theoremstyle{plain}

\numberwithin{equation}{section}

\newtheorem{question}{Question}[section]
\newtheorem{theorem}{Theorem}[section]

\theoremstyle{definition}
\newtheorem{definition}[theorem]{Definition}

\theoremstyle{remark}
\newtheorem{obs}[theorem]{Observation}

\theoremstyle{remark}

\newcounter{commentcounter}

\newcounter{todocounter}

\usepackage{tikz, braids} 
\usepackage[T1]{fontenc} 
\usetikzlibrary{shapes,snakes,positioning,decorations.pathmorphing,calc,backgrounds,chains,fit,matrix,trees,scopes,shapes.geometric,shapes.multipart}
\usetikzlibrary{decorations.pathreplacing, arrows}
\usetikzlibrary{decorations.markings}

\newcommand{\btikz}[1]{\begin{tikzpicture}[#1]}
\newcommand{\etikz}{\end{tikzpicture}}

\pgfdeclarelayer{top}
\pgfsetlayers{main,top}

\newcommand{\poscross}[4]{\begin{scope}[decoration={markings, mark=at position .999 with {\arrow{>}}}]
\draw[thick, looseness=1.5,postaction={decorate}] (1,0) to [out=90, in=-90] (0,1.5);
\draw[white, line width=7, looseness=1.5] (0,0) to [out=90, in=-90] (1,1.5);
\draw[thick, looseness=1.5,postaction={decorate}] (0,0) to [out=90, in=-90] (1,1.5);
\draw (0,0) node[below] {$#1$};
\draw (1,0) node[below] {$#2$};
\draw (0,1.5) node[above] {$#3$};
\draw (1,1.5) node[above] {$#4$};
\end{scope}
}

\newcommand{\negcross}[4]{\begin{scope}[decoration={markings, mark=at position .999 with {\arrow{>}}}]
\draw[thick, looseness=1.5,postaction={decorate}] (0,0) to [out=90, in=-90] (1,1.5);
\draw[white, line width=7, looseness=1.5] (1,0) to [out=90, in=-90] (0,1.5);
\draw[thick, looseness=1.5,postaction={decorate}] (1,0) to [out=90, in=-90] (0,1.5);
\draw (0,0) node[below] {$#1$};
\draw (1,0) node[below] {$#2$};
\draw (0,1.5) node[above] {$#3$};
\draw (1,1.5) node[above] {$#4$};
\end{scope}
}

\newcommand{\twostrandsigmaone}[6]{
\begin{scope}[decoration={markings, mark=at position .2 with {\arrow{>}},mark=at position .8 with {\arrow{>}}}]
\draw[#1,postaction={decorate}] (1,0) node[black,below] {$#4$} to [out=90, in=-90] (2,1.5) node[black,above] {$#4$};

\end{scope}
\begin{scope}[decoration={markings, mark=at position .5 with {\arrow{>}}}]

\draw[#2,postaction={decorate}] (2,0) node[black, below] {$#5$} to [out=90,in=-30] (1.65,.65);
\draw[#3,postaction={decorate}] (1.35,.85) to [out=150,in=-90] (1,1.5) node[black,above] {$#6$};
\end{scope}
}

\newcommand{\twostrandsigmaoneinverse}[6]{
\begin{scope}[yshift=-1.5cm, decoration={markings, mark=at position .5 with {\arrow{>}}}]
\draw[#1,postaction={decorate}] (1,0) to [out=90, in=-150] (1.30,.65);
\draw[#3,postaction={decorate}] (1.7,.85) to [out=30, in=-90] (2,1.5); 
\end{scope}
\begin{scope}[yshift=-1.5cm,decoration={markings, mark=at position .2 with {\arrow{>}}, mark=at position .8 with {\arrow{>}}}]
\draw[#2,looseness=1.5, postaction={decorate}] (2,0) to [out=90,in=-90] (1,1.5);
\end{scope}
\draw (2,-1.5) node[black,below] {$#4$};
\draw (1,0) node[black,above] {$#4$};
\draw (2,0) node[black,above] {$#6$};
\draw (1,-1.5) node[black, below] {$#5$};
}

\author[Delaney]{Colleen Delaney}
\address{Department of Mathematics, UC Santa Barbara, Santa Barbara, CA }
\email{cdelaney@math.ucsb.edu}

\author[Tran]{Alan Tran}
\address{Department of Physics, University of California, Santa Barbara CA}
\email{adtran@physics.ucsb.edu}

\begin{document}

\title{A systematic search of knot and link invariants beyond modular data}

\begin{abstract}

The smallest known example of a family of modular categories that is not determined by its modular data are the rank 49 categories $\mathcal{Z}(\text{Vec}_G^{\omega})$ for $G=\ZZ_{11} \rtimes \ZZ_{5}$. However, these categories can be distinguished with the addition of a matrix of invariants called the $W$-matrix that contains intrinsic information about punctured $S$-matrices. Here we show that it is a common occurrence for knot and link invariants to carry more information than the modular data. We present the results of a systematic investigation of the invariants for small knots and links. We find many small knots and links that are complete invariants of the $\mathcal{Z}(\text{Vec}_G^{\omega})$ when $G=\ZZ_{11} \rtimes \ZZ_{5}$, including the $5_2$ knot.
\end{abstract}

\subjclass[2000]{16W30, 18D10, 19D23}

\date{\today}
\maketitle

\section{Introduction}

Mignard and Schauenburg produced an infinite family of modular tensor categories for which the modular data is not a complete invariant \cite{MS}. This disproved the conjecture that modular tensor categories are determined by their modular data. Motivated by Mignard and Schauenburg, we study invariants of knots and links coming from the smallest counterexample to the conjecture. These are the five inequivalent rank 49 modular tensor categories of the form $\mathcal{Z}(\text{Vec}_G^{\omega})$ when $G$ is the nonabelian group of order 55, $G=\Z_{11}\rtimes Z_{5}$. 

The $S$ and $T$ matrices can be understood as the matrices of invariants associated to admissible labelings of the Hopf link and the twisted unknot, respectively. Thus it is natural to formulate additional categorical invariants from knots and links. In the companion paper \cite{BDGRTW}, we showed that there is an orientation class of the Whitehead link whose associated matrix of invariants, called the $W$-matrix, distinguishes the five categories with the help of the modular data. The motivation for studying the Whitehead link is its relationship to the punctured $S$-matrices, which are a natural extension of the $S$-matrix. A systematic search of link invariants puts this result in context and provides a deeper understanding of the relationship between modular categories and their topological invariants.

We investigate knots up to 10 crossings and links up to 9 crossings that can be represented as the quantum trace of a braidword in $B_3$. It is shown in \cite{BDGRTW} that links coming from $B_2$ are determined by the modular data so we will forego this class of links. We first explore the invariants coming from small knots, highlighting the general behavior with some examples. We show that even the figure-eight knot invariants, taken together with the modular data, is enough to distinguish the categories. On the other hand, the $5_2$ knot distinguishes the categories without the help of the modular data. 

The paper is organized as follows: In section 2 we first discuss what it means for a link to go ``beyond the modular data''. In section 3 we outline the methods of our systematic search and analysis for the rank 49 example. In section 4 we present tables of our results and detail some interesting cases. We conclude in section 5 with some comments and questions about how to extend our results.

This work made use of the Knot Cluster at the Center for Scientific Computing at UCSB, which is supported by NSF Grant CNS-096031. Mathematica packages implementing the calculation of the link invariants and additional data can be found at the authors' websites.
\begin{itemize}
    \item  \url{http://web.physics.ucsb.edu/~adtran/W.html}
\item \url{http://web.math.ucsb.edu/~cdelaney/WMatrices.html}
\end{itemize}

We would like to thank Parsa Bonderson, C\'{e}sar Galindo, Eric Rowell, Zhenghan Wang, Nick Amin, and Sheri Tamagawa for helpful discussions.

\section{Beyond modular data}
Here we make precise what it means for link invariants to go beyond modular data. The discussion applies to any collection of $n$ modular tensor categories which are \emph{not} distinguished by their $S$ and $T$ matrices. That is, there are strictly less than $n$ \emph{sets} of modular data shared among the $n$ modular tensor categories. For a concrete example, one can take the Mignard-Schauenburg categories $\mathcal{Z}(\textrm{Vec}_G^{\omega})$ for $G=\mathds{Z}_{11} \rtimes \mathds{Z}_5$ studied in the next section. We refer the reader to \cite{MS} for more details about the Mignard-Schauenburg categories, and to \cite{BDGRTW} for the general theory of the topological invariants we compute. We also freely use the language of anyon models in our discussion of modular categories, a dictionary for which can be found in \cite{RW18}. Throughout this paper we will assume $G=\mathds{Z}_{11} \rtimes \mathds{Z}_5$ unless otherwise stated.

Let $\{\CC^{(i)}\}$ be such a family of modular categories and fix an ordered basis $\mathcal{L}$ of the anyons. When the modular data are equal as sets, $S^{(i)}$ ($T^{(i)}$) and $S^{(j)}$ ($T^{(j)})$ may not be equal "on the nose", but the entries of the matrices and the multiplicities in which they occur are identical. Moreover, there is a permutation $\rho$ of the anyons so that $S^{(i)}_{\rho(a)\rho(b)}=S^{(j)}_{ab}$ (and $T^{(i)}_{\rho(a)\rho(b)}=T^{(j)}_{ab}$). For if not, one would immediately conclude that there exists no braided tensor auto-equivalence $\tilde{\rho} \in \text{Aut}_{br}^{\otimes}(\CC^{(i)})$ that identifies $\CC^{(i)}$ and $\CC^{(j)}$.

This observation is the key to testing whether a knot or link has the ability to distinguish the categories when taken together with the modular data. Towards this end we define the following notion of a modular permutation

\begin{definition}
    A \emph{modular permutation} $\rho$ between two modular tensor categories $\CC^{(i)}$ and $\CC^{(j)}$ sharing the same sets of modular data is a permutation of the isomorphism classes of simple objects from one category to another
    satisfying
    \begin{align*}
        S^{(i)}_{\rho(a)\rho(b)} &= S^{(j)}_{ab} \\
        T^{(i)}_{\rho(a)\rho(b)} &= T^{(j)}_{ab}
    \end{align*}
    for all $a,b \in \mathcal{L}$. 
\end{definition}
    In particular it follows that the quantum dimensions and fusion rules are left invariant. 
Note that every braided tensor auto-equivalence of a modular category  $\tilde{\rho} \in \text{Aut}_{br}^{\otimes}(\CC)$ induces a modular permutation $\rho$ of its anyons. Of course, not every modular permutation necessarily lifts to a braided tensor auto-equivalence.

Of course one would like to find a link whose associated invariants $L^{(i)}$ distinguish between the $\{\CC^{(i)}\}$ as sets. This is a stronger notion of what it means to go "beyond modular data". But observe that it is not necessary for a link to have this property for it to distinguish between a collection of modular tensor categories when the modular data is known. In fact, a link invariant could be identical "on the nose" for all $\CC^{(i)}$, and yet together with the $S$ and $T$ matrices distinguish between members of the collection: if there is no modular permutation that extends to the $L$ invariants, one can immediately conclude that none of the $\CC^{(i)}$ are related by a braided-tensor autoequivalence. The figure eight knot has this property for the rank 49 Mignard-Schauenburg categories, and we describe the phenomenon in detail in section 3.1.

In light of these considerations, we define two notions of distinguishability.

\begin{definition}
 A framed link $L$ distinguishes a family of modular categories $\{\CC^{(i)}\}$ 

\begin{enumerate}

\item \emph{weakly} if there exists $i,j$ such that none of the modular permutations $\{\rho\}^{(i,j)}$ map $L^{(i)}$ onto $L^{(j)}$,

\item \emph{strongly} if $L^{(i)} \ne L^{(j)}$ as sets for all $i\ne j$.

\end{enumerate}
\end{definition}
Note that even if $L$ weakly distinguishes there may still be strictly less than $n$ distinct sets among $\{S^{(i)}, T^{(i)}, L^{(i)}\}$; it is simply the statement that permutations between $\CC^{(i)}$ and $\CC^{(j)}$ are not simultaneously compatible for $S$, $T$, and $L$. It is clear that strongly distinguishing implies weakly distinguishing. 

Our main result is the classification of which knots and two-component links with braid word representatives in $B_3$ either weakly or strongly distinguish $\mathcal{Z}(\text{Vec}_G^{\omega})$ for $G=\ZZ_{11} \rtimes \ZZ_{5}$.
\begin{theorem}
The rank 49 Mignard-Schauenburg categories are 

strongly distinguished by

\begin{itemize}
    \item the knots $5_2,8_{n21},10_2,10_{46},10_{94},10_{106}, 10_{n126}, 10_{n155}$ and
    \item the links $6^2_3,7^2_1,7^2_2+-,7^2_5+-,8^2_{11},8^2_3,9^2_2,9^2_{20},9^2_{23},9^2_{21}+-,9^2_{34},9^2_{39},9^2_{51},9^2_{52},9^2_{54}+-,9^2_{58}+-,9^2_{59}+-$
    \end{itemize}

and weakly distinguished by

\begin{itemize}
    \item the knots $4_1, 8_9,8_{18}$ and
    \item the links $5_1^2,7^2_4,7^2_6,7^2_8, 9^2_5, 9^2_{13}, 9^2_{31},9^2_{37},9^2_{41},9^2_{44},9^2_{50},9^2_{55},9^2_{57}+-$.
    \end{itemize}
\end{theorem}




\subsection{$B$-type anyons and quandle coloring numbers}
When $G=\mathds{Z}_{11} \rtimes \mathds{Z}_5$, there are three types of anyons in $\mathcal{Z}(\text{Vec}_G^{\omega})$, which we call $I$-type, $A$-type, and $B$-type. In \cite{BDGRTW} it was shown that the invariants coming from labeling all components of a link by a fixed type $B$-anyon can be realized as certain quandle coloring numbers. We will see that the $B$-type anyons play a special role in the Mignard-Schauenburg categories, and that it is useful to be able to compute their invariants in more than one way.

We recall below the theorem that relates the topological invariants to certain quandle coloring numbers.  

The quandle in question is given as follows. Given $G=\Z_q\rtimes _n\Z_p=\langle a\rangle\rtimes_n \langle b\rangle$, the quandles associated to the conjugacy classes of $b^k,k=0,1,...,p-1$ are Alexander quandles: 
\begin{align*}
    X_k=\Z/q\Z,&&  a\triangleright b:=(1-n^k)a+n^kb 
\end{align*}
for  $a, b \in \Z/q\Z$.

Here $b$ is a generator of the copy of $\mathds{Z}_5$ in $G=\mathds{Z}_{11} \rtimes \mathds{Z}_5$. The parametrization $([b^k],\pi_k^s)$ corresponds to a type $B_{k,s}$ anyon.
\begin{theorem}
Let $\sigma_{i_1}^{\epsilon_1}\sigma_{i_2}^{\epsilon_2}\cdots \sigma_{i_h}^{\epsilon_h} \in \cB_n$, and $L=\widehat{\sigma}$ its closure. The invariant of the oriented framed link $L$ colored by the simple object $V(b^k,\pi_k^s)$ is
\begin{equation}
    q^{Wr(L)}C_{X_k}(L),
\end{equation}
where $Wr(L)$ is the writhe of $L$, $X_k$ is the Alexander quandle constructed above and $q=e^{\frac{2\pi i}{p^2}(sp+uk )k}$, which is the twist $\theta_{([b^k],\pi_k^s)}$ of the anyon $([b^k],\pi_k^s)$.
\end{theorem}
In particular, when the write is zero, the invariant is just an integer.

\section{Methods of systematic search}
We calculated the invariants for oriented knots up to 10 crossings and links up to 9 crossings that can be realized with three strands or less. Calculating the invariants for $B_n$ where $n>3$ is possible, but the mathematical theory is more complicated to implement.

So-called minimum braid word representatives of knots and links for small crossing numbers were tabulated by Gittings in \cite{gittings}. A minimum braid word representative $b\in B_n$ of a link $L$ uses the smallest possible number of strands and the shortest possible length as a word in the braid group generators $\sigma_i$. The braids in Gitting's paper are required to satisfy other conditions so that they are uniquely defined, but for our purposes it suffices to use any braid word $b$ that minimizes $n$ and $|b|$. This enables one to compute the invariants in an efficient manner. Gitting's table provides a minimum braid word from each orientation class of the knots and links up to orientation reversal. 

Since we are interested in which invariants beyond modular data, it is not necessary to analyze those knots and links that can be realized as trace closures of two-strand braids. It was shown in \cite{BDGRTW} that closures of two-strand braids can be expressed entirely in terms of the modular data. Therefore any symmetry of the modular data is a symmetry of a two-strand invariant. 

This leaves the three-strand braid words that correspond to oriented links. Gittings' convention for labeling braid word representatives in $B_3$ uses  $A$ and $B$ to denote generators of the braid group, with $a$ and $b$ their respective inverses. Braid diagrams are read from the top down. However, these differ slightly from the typical conventions when studying braid group representations from modular categories, where the right-handed crossings are the generators and the left-handed crossings their inverses. Moreover, we read braid diagrams from the bottom up.

$$A= \begin{tikzpicture}[baseline=.5cm]
\negcross{}{}{}{}
\begin{scope}[thick, decoration={markings, mark=at position .999 with {\arrow{>}}}]
\draw[postaction={decorate}] (2,0)--(2,1.5);
\end{scope}
\end{tikzpicture} \hspace{25pt},\hspace{25pt} B= \begin{tikzpicture}[baseline=.5cm]
\negcross{}{}{}{}
\begin{scope}[thick, decoration={markings, mark=at position .999 with {\arrow{>}}}]
\draw[postaction={decorate}] (-1,0)--(-1,1.5);
\end{scope}
\end{tikzpicture}
$$

For example, the braid $AAABaB$ corresponds to 
$\sigma_1^{-3} \sigma_2^{-1} \sigma_1 \sigma_2^{-1}$ which has braid closure
$$
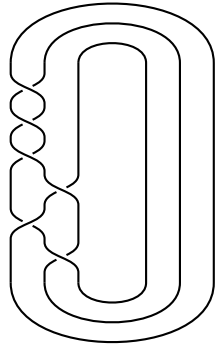
\begin{tikzpicture}[scale=.45]
\begin{scope}[thick]
\braid a_1 a_1 a_1 a_2 a_1^{-1} a_2 ;

\draw (3,-6.5) to [out=-90, in=-90] (5,-6.5)--(5,0) to  [out=90,in=90] (3,0);
\draw (2,-6.5) to [out=-90, in=-90]  (6,-6.5)--(6,0) to  [out=90,in=90] (2,0);
\draw (1,-6.5) to [out=-90, in=-90] (7,-6.5)--(7,0) to  [out=90,in=90] (1,0);
\end{scope}
\end{tikzpicture}.$$

\subsection{Modular permutations}
The set of modular permutations between $\mathcal C_i$ and $\mathcal C_j$ can be found in a straight-forward way. By matching up quantum dimensions, twists and fusion rules.
With this restriction one can exhaustively search through all $T$-invariant permutations. The ones which also apply to the $S$-matrices are the set of modular permutations.

For example, in the MS example there are five categories indexed by  $\mathcal{Z}(\text{Vec}_G^{\omega^u})_{0\leq u \leq 4}$. The $u=0$, $u=1,4$ and $u=2,3$ categories share the same modular data sets.

There are 49 anyons, which we organized into three types: type I-anyons, type A-anyons, and type B-anyon. There are seven type I-anyons, with $\{I_0,I_1, \ldots, I_4\}$ having quantum dimension 1 and $\{I_5,I_6\}$ having quantum dimension 5. There are twenty-two type $A$ anyons which split into two classes of eleven each, labeled $\{A_{1,i},A_{2,j}\}$ for $0\le i,j \le 10$. The type $A$ anyons all have quantum dimension 5. There are twenty type $B$ anyons, that split into four blocks of five each, labeled $\{B_{1,i},B_{2,j}, B_{3,k},B_{4,l}\}$, where $0 \le  i,j,k,l \le 4$. Each $B$-type anyon has quantum dimension 11. We refer the reader to \cite{MS,BDGRTW} for more details about the category.

Thus five anyons have quantum dimension 1, twenty-four have dimension 5 and twenty have dimension 11. Further, one can verify that each non-trivial twist occurs exactly twice among the dimension 5 and 11 anyon types. 

Since the identity object must always be permuted to itself there are $4!$ ways to permute the remaining dimension 1 objects. Of the twenty-four dimension 5 anyons, 4 have trivial twist and so there are an additional $4!$ permutations of these. The twenty leftover dimension 5 anyons are paired up according to their twists. Each pair can either swap their members or not under the action of a permutation, so there are $2^{10}$ choices here. 
Finally, the twenty dimension 11 anyons are naturally split into four blocks of five according to their fusion rules. The blocks further pair off according to their twists. The pairs of blocks can now either simultaneously swap all their elements, or not; so there are $2^2$ permutations here. If the permutations do not conform to this block structure the fusion structure will not be preserved.

There is thus a total of $4!4!2^{10}2^2=2359296$ $T$-respecting permutations. Only 8 of them lift to modular permutations between $u=1,4$ and likewise between $u=2,3$. Tables listing the modular permutations are given in Appendix B. 


Having fixed the basis $$\{\{I_r\},\{A_{1,i}\},\{A_{2,i}\}, \{B_{1,k}\},\{B_{2,k}\}, \{B_{3,k}\}, \{B_{4,k}\}\}_{{0\leq r \leq 6 \atop 0\leq i \leq10} \atop 0\leq k\leq 4},$$ we represent the invariants of a knot with a vector of 49 entries, corresponding to each of the admissible labelings by anyons. For two-component links we express the invariants as a $49 \times 49$ matrix written in the given basis. For three-component links, of which we only consider one example, the invariants are organized into a $49 \times 49 \times 49$ tensor.
We generate the invariant data for each knot and link and classify whether it weakly or strongly distinguishes the five categories of $\mathcal{Z}(\text{Vec}_G^{\omega^u})$.


\section{Results of systematic search}
We provide tables of the results, first for knots and then for two-component links. 

For knots, the first column contains the Rolfsen ID \cite{Rolfsen} of the link whose representative braidword is in the second column. Then we indicate whether the knot distinguishes weakly or strongly. While strongly distinguishing implies weakly distinguishing, we only mark the final column when a knot or link strongly distinguishes.

\subsection{Knots}
\begin{center}
\begin{align}
\hspace{-2.6em}
\scalemath{0.9}{
\begin{array}{c|c|c|c}
 \text{Knot} & \text{Braidword} & \text{Weakly} & \text{Strongly} \\
\hline \hline
4_1 & \text{AbAb} & \checkmark  &   \\
\hline
 5_2 & \text{AAABaB} &   & \checkmark  \\
\hline
 6_2 & \text{AAAbAb} &   &   \\
\hline
 6_3 & \text{AAbAbb} &   &   \\
\hline
 7_3 & \text{AAAAABaB} &   &   \\
\hline
 7_5 & \text{AAAABaBB} &   &   \\
\hline
 8_2 & \text{AAAAAbAb} &   &   \\
\hline
 8_7 & \text{AAAAbAbb} &   &   \\
\hline
 8_5 & \text{AAAbAAAb} &   &   \\
\hline
 8_{10} & \text{AAAbAAbb} &   &   \\
\hline
 8_9 & \text{AAAbAbbb} & \checkmark  &   \\
\hline
 8_{16} & \text{AAbAAbAb} &   &   \\
\hline
 8_{17} & \text{AAbAbAbb} &   &   \\
\hline
 8_{18} & \text{AbAbAbAb} & \checkmark  &   \\
\hline
 8_{n20} & \text{AAAbaaab} &   &   \\
\hline
 8_{n19} & \text{AAABAAAB} &   &   \\
\hline
 8_{n21} & \text{AAABaaBB} &   & \checkmark  \\
\hline
 9_{29} & \text{AAAAAAABaB} &   &   \\
\hline
9_6 & \text{AAAAAABaBB} &   &   \\
\hline
9_9 & \text{AAAAABaBBB} &   &   \\
\hline
 9_{16} & \text{AAAABBaBBB} &   &   \\
\hline
10_2 & \text{AAAAAAAbAb} &   & \checkmark  \\
\hline
 10_5 & \text{AAAAAAbAbb} &   &   \\
\hline
 10_{46} & \text{AAAAAbAAAb} &   & \checkmark  \\
\hline
 10_{47} & \text{AAAAAbAAbb} &   &   \\
\hline
10_9 & \text{AAAAAbAbbb} &   &   \\
\hline
 10_{62} & \text{AAAAbAAAbb} &   &   \\
\hline
 10_{85} & \text{AAAAbAAbAb} &   &   \\
 \hline
10_{82} & \text{AAAAbAbAbb} &   &   \\
\end{array}
\qquad \qquad
\begin{array}{c|c|c|c}
 \text{Knot} & \text{Braidword} & \text{Weakly} & \text{Strongly} \\
\hline \hline
 10_{17} & \text{AAAAbAbbbb} &   &   \\
\hline
 10_{48} & \text{AAAAbbAbbb} &   &   \\
\hline
 10_{64} & \text{AAAbAAAbbb} &   &   \\
\hline
 10_{100} & \text{AAAbAAbAAb} &   &   \\
\hline
 10_{94} & \text{AAAbAAbbAb} &   & \checkmark  \\
\hline
 10_{106} & \text{AAAbAbAAbb} &   & \checkmark  \\
\hline
 10_{112} & \text{AAAbAbAbAb} &   &   \\
\hline
10_{91} & \text{AAAbAbbAbb} &   &   \\
\hline
 10_{79} & \text{AAAbbAAbbb} &   &   \\
\hline
10_{104} & \text{AAAbbAbAbb} &   &   \\
\hline
 10_{116} & \text{AAbAAbAbAb} &   &   \\
\hline
 10_{99} & \text{AAbAAbbAbb} &   &   \\
\hline
10_{118}& \text{AAbAbAbbAb} &   &   \\
\hline
 10_{109} & \text{AAbAbbAAbb} &   &   \\
\hline
 10_{123} & \text{AbAbAbAbAb} &   &   \\
\hline
 10_{n125}& \text{AAAAAbaaab} &   &   \\
\hline
 10_{n124} & \text{AAAAABAAAB} &   &   \\
\hline
 10_{n126} & \text{AAAAABaaaB} &   & \checkmark  \\
\hline
10_{n127} & \text{AAAAABaaBB} &   &   \\
\hline
10_{n141}& \text{AAAAbaaabb} &   &   \\
\hline
 10_{n139} & \text{AAAABAAABB} &   &   \\
\hline
 10_{n143} & \text{AAAABaaaBB} &   &   \\
\hline
 10_{n148} & \text{AAAABaaBaB} &   &   \\
\hline
 10_{n149} & \text{AAAABaBaBB} &   &   \\
\hline
 10_{n155} & \text{AAABaaBaaB} &   & \checkmark  \\
\hline
 10_{n161} & \text{AAABaBAABB} &   &   \\
\hline
 10_{n159} & \text{AAABaBaaBB} &   &   \\
\hline
10_{n152} & \text{AAABBAABBB} &   &   \\
\hline
 10_{n157} & \text{AAABBaBaBB} &   &   \\
\end{array}}
\notag
\end{align}
\captionof{table}{List of knots up to 10 crossings with representative in $B_3$. Whether the knot \emph{weakly} or \emph{strongly} distinguishes $\mathcal{Z}(\text{Vec}_G^{\omega^u})$ is indicated with $\checkmark$.}
\end{center}

The general behavior of small 3-strand knot invariants can be neatly summarized according to the following observations.
\begin{obs}
\hspace{1em}
\begin{enumerate}
\item  The knot invariants of $I$-type anyons are all integral and independent of the 3-cocycle $\omega$.

\item The knot invariants of $A$-type anyons are independent of the 3-cocycle $\omega$.
\item  The invariants of knots with zero writhe are all equal and independent of the 3-cocycle $\omega$. 

\end{enumerate}
\end{obs}

Thus any hope for a knot to strongly distinguish the five categories lies with the type $B$-anyons. However, as the table above shows that it is not too rare for a knot invariant to distinguish the Mignard-Schauenburg categories.

Next we treat the details of some specific knots that exhibit the general phenomena found in the search. 
\subsubsection{The figure eight knot goes beyond modular data}

A direct calculation of the invariants assigned to the figure eight knot $4_1$ shows that they are identical for all 5 choices of cocycle $\omega$. However, the vector of invariants fails to be left invariant under permutations of the modular data, so while it cannot detect the cocycle $\omega$, even the figure eight knot invariant goes ``beyond the modular data".

To see this consider the vector of invariants colored by $B$-type anyons.
\begin{equation*}
\begin{array}{c|c|c}

\text{Anyon type} & \text{Conj. Class} & \text{Number of admissible colorings} \\
\hline



B_1 & [b^1] & 11 \\

B_2 & [b^2] & 11 \\

B_3 & [b^3] & 11 \\

B_4 & [b^4] & 11

\end{array}
\end{equation*}
where $B_{k}$ is the set of $\{B_{k,s}\}_{0\leq s\leq 4}$ anyons.

By inspection of the twists (see Appendix), any equivalence of the categories for  $u=1$ and $u=4$ would have to send 
\begin{equation*}
    \begin{array}{ccc}
    B^{(1)}_1 &\mapsto & B^{(4)}_3 \\
    B^{(1)}_2 &\mapsto & B^{(4)}_1 \\
    B^{(1)}_3 &\mapsto & B^{(4)}_4 \\
    B^{(1)}_4 &\mapsto & B^{(4)}_2 
    \end{array}
    \qquad \text{or} \qquad
    \begin{array}{ccc}
    B^{(1)}_1 &\mapsto & B^{(4)}_2 \\
    B^{(1)}_2 &\mapsto & B^{(4)}_4 \\
    B^{(1)}_3 &\mapsto & B^{(4)}_1 \\
    B^{(1)}_4 &\mapsto & B^{(4)}_3 
    \end{array}
\end{equation*}

But the integer-valued invariants of the figure eight vector are different for each pairing. Therefore the two categories for $u=1$ and $u=4$ are inequivalent. A similar result holds for $u=2$ and $u=3$.

We remark that among braids requiring at least three strands the figure eight knot has the braid word representative of shortest length. So even the ``smallest" knot is beyond modular data.

As a further check on the $B$-labeled invariants, we can take advantage of the realization of $B$-type invariants as a certain quandle coloring number. \cite{BDGRTW}
We compute the invariants $C_{X_k}$ for the figure-eight knot explicitly for the four quandles $(X_k,*)$ where $X_k=[b^k]$ for $1\le k \le 4$ with $*$ given by $x * y = yxy^{-1}$ and $\bar{*}$ given by $x \bar{*} y = y^{-1}xy$. 

Note that here we read braid words bottom to top, and use the convention that
$$\sigma_1= \begin{tikzpicture}[baseline=.5cm]
\poscross{}{}{}{}
\end{tikzpicture}\hspace{25pt},\hspace{25pt} \sigma_1^{-1}= \begin{tikzpicture}[baseline=.5cm]
\negcross{}{}{}{}
\end{tikzpicture}
$$
The underlying sets of the quandles are given by the conjugacy classes $[b^k]$, corresponding to the four classes of $B$-anyons.

$$\begin{array}{c|c}

\textrm{Anyon type} & \textrm{Conjugacy class }[b^k] \subset G \\
\hline

B_1  & \{a^lb\}_{l=0}^{10} \\
B_2  & \{a^lb^2\}_{l=0}^{10} \\
B_3  & \{a^lb^3\}_{l=0}^{10} \\
B_4  & \{a^lb^4\}_{l=0}^{10} \\
\end{array}$$

We calculate the admissible quandle colorings of the link $\widehat{(\sigma_1^{-1} \sigma_2)^2}$.

There are four arcs in this knot. We start by labeling the four arcs, $x$ (red), $y$ (blue), $z$ (black), $w$ (green). Then order of the labeling proceeds from right to left, since the first generator is $\sigma_2$. 

There are four crossings $$
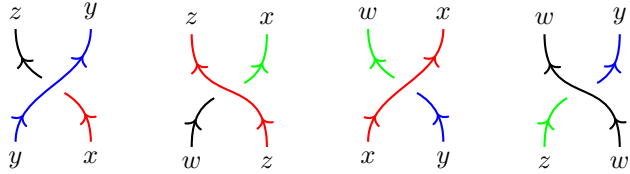
\begin{tikzpicture}[thick]\twostrandsigmaone{blue}{red}{black}{y}{x}{z}\end{tikzpicture}\hspace{25pt}
\begin{tikzpicture}[thick] \twostrandsigmaoneinverse{black}{red}{green}{z}{w}{x}
\end{tikzpicture}
\hspace{25pt}
\begin{tikzpicture}[thick]\twostrandsigmaone{red}{blue}{green}{x}{y}{w}\end{tikzpicture}\hspace{25pt}\btikz{thick}
\twostrandsigmaoneinverse{green}{black}{blue}{w}{z}{y}
\etikz$$
corresponding to the relations $z=x * y$, $x=w \bar{*} z$, $w=y * x $, and $y=z\bar{*} w$, respectively. It follows from the first and third constraints that $z$ and $w$ are determined by $x$ and $y$. 
Then admissible colorings of the arcs are in correspondence with pairs $(x,y)$ where $x,y \in [b^k]$ and the equations

$$x=(y * x) \bar{*} (x * y)$$
and 
$$y= (x * y) \bar{*} (y * x)$$
are satisfied. 

These equations in $G$ were solved by computer, with the following results. When $k=1$ or $k=4$, these equations only hold in $G=\ZZ_{11}\rtimes Z_5$ when $x=y$. In this case there are $|[b^k]|=11$ admissible colorings of the link by the quandle. 
When $k=2$ or $k=3$, these equations are satisfied by every $x,y \in [b^k]$ and hence there are $|[b^k]|^2=121$ admissible colorings. This corresponds exactly to the calculation of the invariants via the representations.



Observe that these relations have an $x-y$ symmetry: the two relations can be obtained from one another by exchanging $x$ and $y$.














\subsection{Two-component links}
The data is organized similarly to table 1, except now the first column contains the Doll-Hoste ID \cite{DH} of the link with representative braid word in the second column.

\begin{center}
\begin{align*}
\hspace{-2.6em}
\scalemath{0.9}{
\begin{array}{l|l|c|c}
 \text{Link} & \text{Braidword} & \text{Weakly} & \text{Strongly} \\
\hline\hline
4^2_1 & \text{AABaB} &   &   \\
\hline
5^2_1 & \text{AAbAb} & \checkmark  &   \\
\hline
 6^2_3 & \text{AAABaBB} & &\checkmark \\
\hline
 6^2_2 & \text{AAAABaB} &   &   \\
\hline
 7^2_1 & \text{AAAAbAb} & &\checkmark \\
\hline
 7^2_4 & \text{AAAbAAb} & \checkmark  &   \\
\hline
 7^2_2+1 & \text{AAAbAbb} & &\checkmark \\
\hline
 7^2_5+- & \text{AAbAAbb} & &\checkmark \\
\hline
7^2_6 & \text{AAbAbAb} & \checkmark  &   \\
\hline
7^2_7 & \text{AAAbaab} &   &   \\
\hline
 7^2_7+- & \text{AAABAAB} &   &   \\
\hline
 7^2_8 & \text{AAABaaB} & \checkmark  &   \\
\hline
8^2_{11} & \text{AAAABBaBB} & &\checkmark \\
\hline
8^2_2 & \text{AAAAAABaB} &   &   \\
\hline
 8^2_3 & \text{AAAAABaBB} & &\checkmark \\
\hline
 8^2_4+- & \text{AAAABaBBB} &   &   \\
\hline
9^2_1 & \text{AAAAAAbAb} &   &   \\
\hline
 9^2_{13} & \text{AAAAAbAAb} & \checkmark  &   \\
\hline
 9^2_2 & \text{AAAAAbAbb} & &\checkmark \\
\hline
 9^2_{19} & \text{AAAAbAAAb} &   &   \\
\hline
 9^2_{20} & \text{AAAAbAAbb} & &\checkmark \\
\hline
 9^2_{31} & \text{AAAAbAbAb} & \checkmark  &   \\
\hline
 9^2_5 & \text{AAAAbAbbb} & \checkmark  &   \\
\hline
 9^2_{14}+- & \text{AAAAbbAbb} &   &   \\
\end{array}
\qquad\qquad
\begin{array}{l|l|c|c}
 \text{Link} & \text{Braidword} & \text{Weakly} & \text{Strongly} \\
\hline\hline
 9^2_{23} & \text{AAAbAAAbb} & &\checkmark \\
\hline
 9^2_{35} & \text{AAAbAAbAb} &   &   \\
\hline
 9^2_{21}+- & \text{AAAbAAbbb} & &\checkmark \\
\hline
 9^2_{34} & \text{AAAbAbAbb} & &\checkmark \\
\hline
 9^2_{37} & \text{AAAbAbbAb} & \checkmark  &   \\
\hline
 9^2_{29}+- & \text{AAAbbAAbb} &   &   \\
\hline
 9^2_{40}+- & \text{AAbAAbAAb} &   &   \\
\hline
 9^2_{39} & \text{AAbAAbAbb} & &\checkmark \\
\hline
 9^2_{41} & \text{AAbAbAAbb} & \checkmark  &   \\
\hline
 9^2_{42} & \text{AAbAbAbAb} &   &   \\
\hline
9^2_{43} & \text{AAAAAbaab} &   &   \\
\hline
 9^2_{43}+- & \text{AAAAABAAB} &   &   \\
\hline
 9^2_{44} & \text{AAAAABaaB} & \checkmark  &   \\
\hline
9^2_{50} & \text{AAAAbaaab} & \checkmark  &   \\
\hline
 9^2_{49} & \text{AAAABAAAB} &   &   \\
\hline
 9^2_{51} & \text{AAAABaaaB} & &\checkmark \\
\hline
 9^2_{52} & \text{AAAABaaBB} & &\checkmark \\
\hline
9^2_{55} & \text{AAAABaBaB} & \checkmark  &   \\
\hline
 9^2_{54}+- & \text{AAAbaaabb} & &\checkmark \\
\hline
 9^2_{53} & \text{AAABAAABB} &   &   \\
\hline
 9^2_{57}+- & \text{AAABaaBaB} & \checkmark  &   \\
\hline
 9^2_{58}+- & \text{AAABaBaBB} & &\checkmark \\
\hline
 9^2_{59}+- & \text{AAABBAABB} & &\checkmark \\
\hline
 9^2_{61} & \text{AABaBAABB} &   &  
\end{array}
}
\end{align*}
\captionof{table}{List of two-component links up to 9 crossings with representative in $B_3$. Whether the link \emph{weakly} or \emph{strongly} distinguishes $\mathcal{Z}(\text{Vec}_G^{\omega^u})$ is indicated with $\checkmark$.}
\end{center}

\subsubsection{The Whitehead link}
The invariants associated to the Whitehead link were the focus of \cite{BDGRTW}, where the authors defined the $W$-matrix related to the Whitehead invariants by 

$$
W_{ab}=\frac{\theta_a}{\theta_b}\widetilde{W}_{ab}=\frac{\theta_a}{\theta_b}
\begin{tikzpicture}[baseline=0, thick,scale=.33, shift={(0,-4.8)}]
%
%
\begin{scope}[decoration={markings, mark=at position 0.5 with {\arrow{>}}}]
\draw [looseness=1](2.,6.75) to [out=180, in=90] (1.25,4.75);
\draw [looseness=1,white,line width=5](1.5,6) to [out=180, in=-90] (0.75,7.5);
\draw [looseness=1](1.5,6) to [out=180, in=-90] (0.75,7.5);
\draw [looseness=1](1.5,6) to [out=0, in=-90] (2.25,7.5);
\draw [looseness=1,white,line width=5](2,6.75) to [out=0, in=90 ] (2.75,4.75);
\draw [looseness=1](2.,6.75) to [out=0, in=90] (2.75,4.75);
\draw (2.75,4.75) -- (2.75,2.0);
\draw (1.25,4.75) -- (1.25,2.0);
\draw (2.75,2.0) to [out=-90, in=-90] (5,2.0);
\draw (1.25,2.0) to [out=-90, in=-90] (6.5,2.0);
\draw[postaction=decorate] (0.75,7.5) -- (0.75,8);
\draw (0.75,8) node[left] {\small $a$};
\draw (2.25,7.5) -- (2.25,8);
\draw (2.25,8) to [out=90, in=90] (5,8);
\draw (0.75,8) to [out=90, in=90] (6.5,8);
\draw (5,2) -- (5,8);
\draw (6.5,2) -- (6.5,8);
\end{scope}
\begin{scope}[decoration={markings, mark=at position 0.56 with {\arrow{<}}}]
\draw[white, line width=5] (2,4) ellipse (1.75cm and 0.5cm);
\draw[postaction={decorate}] (2,4) ellipse (1.75cm and 0.5cm);
\draw (0.5,3.6) node[below] {\small $b$};
\draw[white, line width=5] (1.25,5) -- (1.25,4);
\draw[white, line width=5] (2.75,5) -- (2.75,4);
\draw (1.25,5) -- (1.25,4);
\draw (2.75,5) -- (2.75,4);
\end{scope}
\end{tikzpicture}\,\,.
%
$$

The $W$-matrix contains intrinsic information about the punctured $S$-matrices of a modular category, which makes it a natural candidate for an invariant to complement the modular data. A detailed proof of the weakly distinguishing property of the Whitehead link was given in \cite{BDGRTW}.

However, our results suggest that it is common for three-strand two-component links for $\mathcal{Z}(\text{Vec}_G^{\omega^u})$ will be beyond the modular data. Moreover, several of the two-component links \emph{strongly} distinguish the family of modular categories.

\subsection{Three-component links}
Due to computational complexity, we focus here on a single three-component link, the Borromean rings.

\subsubsection{The Borromean rings}

The invariants of the Borromean rings $6^3_2$ behave identically to the invariants of the figure eight knot $4_1$: even though the invariants are independent of the cocycle $\omega$, they still weakly distinguish the Mignard-Schauenburg categories. During the preparation of this manuscript we learned that the weakly-distinguishing property of the Borromean rings was independently discovered by Schauenburg and collaborators \cite{PSPrivate}.

\begin{tabular}{l|l|c|c|c}

\text{Link} & \text{Braidword} $b$ & \text{All equal} & \text{Weakly distinguishes}  & \text{Strongly distinguishes} \\

\hline
\hline
$4_1$ & AbAb & $\checkmark$ & $\checkmark$ &  \\
\hline
$6^3_2$ & AbAbAb & $\checkmark$& $\checkmark$ &  \\

\end{tabular}

Since the writhe of a braid representative of the Borromean rings is zero, Theorem 2.2 ensures that the invariants associated to labeling all components by $B$-type anyons are all integral. Indeed, we find the invariants are all equal to 11, both by direct calculation and by counting quandle colorings. 

\begin{equation*}
\begin{array}{c|c|c}

\text{Anyon type} & \text{Conj. Class} & \text{Number of admissible colorings} \\
\hline



B_1 & [b^1] & 11 \\

B_2 & [b^2] & 11 \\

B_3 & [b^3] & 11 \\

B_4 & [b^4] & 11

\end{array}
\end{equation*}

\section{Discussion of results and open questions}

\subsection{The Mignard-Schauenburg categories $\mathcal{Z}(\text{Vec}_{\ZZ_{q} \rtimes \ZZ_p}^{\omega}))$}

It has been checked by computer that all modular categories up to rank 32 are distinguished by their modular data \cite{MS}. We have shown that for the smallest known modular category for which the modular data is not complete, there are many knots and links of few crossings which distinguish them. However, Mignard and Schauenburg produced an infinite family of counterexamples. It is not known whether the same set of knots and links distinguish the larger counterexamples as well, but for the next few members in the family it is something that could be checked with our methods. 

\subsection{Invariants of general modular categories}
Now that we know that small knots and links can contain powerful information about large categories, it is natural to ask the following question.

\begin{question}
Does there exists a set of framed links $\{L_k\}$ such that $\{L_k\}$ is a complete invariant for all modular tensor categories?
\end{question}

This question is somewhat out of reach with computational methods, as calculation of link invariants for a modular tensor category $\CC$ will in general require the full algebraic data $\{N_a^{bc}, R_a^{bc}, [F^{abc}_d]_{n;m}\}$. Solving the pentagon and hexagon equations that must be satisfied by the braiding and associators of $\CC$ is a notoriously hard problem, even for small categories and with the help of a computer. In our case, we were able to calculate the invariants because they arose as traces of certain representations that can be calculated using only the data of the group $G$ and its third cohomology $H^3(G;U(1))$.

\appendix
\section{$T$-matrices}
\begin{equation}\notag
\scalemath{0.91}{
	\begin{array}{c|cc}
	\textrm{Label} & d & \theta \\ \hline \hline
	I_0 & 1 & 1\\
	I_1 & 1 & 1\\
	I_2 & 1 & 1\\
	I_3 & 1 & 1\\
	I_4 & 1 & 1\\
	I_5 & 5 & 1\\
	I_6 & 5 & 1\\ \hline
	A_{1,0} & 5 & 1\\
	A_{1,1} & 5 & \exp(\frac{i2\pi}{11})\\
	A_{1,2} & 5 & \exp(\frac{i2\pi}{11}2)\\
	A_{1,3} & 5 & \exp(\frac{i2\pi}{11}3)\\
	A_{1,4} & 5 & \exp(\frac{i2\pi}{11}4)\\
	A_{1,5} & 5 & \exp(\frac{i2\pi}{11}5)\\
	A_{1,6} & 5 & \exp(\frac{i2\pi}{11}6)\\
	A_{1,7} & 5 & \exp(\frac{i2\pi}{11}7)\\
	A_{1,8} & 5 & \exp(\frac{i2\pi}{11}8)\\
	A_{1,9} & 5 & \exp(\frac{i2\pi}{11}9)\\
	A_{1,10} & 5 & \exp(\frac{i2\pi}{11}10)\\ \hline
	A_{2,0} & 5 & 1\\
	A_{2,1} & 5 & \exp(\frac{i4\pi}{11})\\
	A_{2,2} & 5 & \exp(\frac{i4\pi}{11}2)\\
	A_{2,3} & 5 & \exp(\frac{i4\pi}{11}3)\\
	A_{2,4} & 5 & \exp(\frac{i4\pi}{11}4)\\
	A_{2,5} & 5 & \exp(\frac{i4\pi}{11}5)\\
	A_{2,6} & 5 & \exp(\frac{i4\pi}{11}6)\\
	A_{2,7} & 5 & \exp(\frac{i4\pi}{11}7)\\
	A_{2,8} & 5 & \exp(\frac{i4\pi}{11}8)\\
	A_{2,9} & 5 & \exp(\frac{i4\pi}{11}9)\\
	A_{2,10} & 5 & \exp(\frac{i4\pi}{11}10)\\ \hline
	B_{1,0} & 11 & \exp(\frac{i2\pi}{25} 1^2 u) \\
	B_{1,1} & 11 & \exp(\frac{i2\pi}{25}(5\cdot 1\cdot 1+1^2 u)) \\
	B_{1,2} & 11 & \exp(\frac{i2\pi}{25}(5\cdot 1\cdot 2+1^2 u)) \\
	B_{1,3} & 11 & \exp(\frac{i2\pi}{25}(5\cdot 1\cdot3+1^2 u)) \\
	B_{1,4} & 11 & \exp(\frac{i2\pi}{25}(5\cdot 1\cdot4+1^2 u)) \\ \hline
	B_{2,0} & 11 & \exp(\frac{i2\pi}{25} 2^2 u) \\
	B_{2,1} & 11 & \exp(\frac{i2\pi}{25}(5\cdot 2\cdot1+2^2 u)) \\
	B_{2,2} & 11 & \exp(\frac{i2\pi}{25}(5\cdot 2\cdot2+2^2 u)) \\
	B_{2,3} & 11 & \exp(\frac{i2\pi}{25}(5\cdot 2\cdot3+2^2 u)) \\
	B_{2,4} & 11 & \exp(\frac{i2\pi}{25}(5\cdot 2\cdot4+2^2 u)) \\ \hline
	B_{3,0} & 11 & \exp(\frac{i2\pi}{25} 3^2 u) \\
	B_{3,1} & 11 & \exp(\frac{i2\pi}{25}(5\cdot 3\cdot1+3^2 u)) \\
	B_{3,2} & 11 & \exp(\frac{i2\pi}{25}(5\cdot 3\cdot2+3^2 u)) \\
	B_{3,3} & 11 & \exp(\frac{i2\pi}{25}(5\cdot 3\cdot3+3^2 u)) \\
	B_{3,4} & 11 & \exp(\frac{i2\pi}{25}(5\cdot 3\cdot4+3^2 u)) \\ \hline
	B_{4,0} & 11 & \exp(\frac{i2\pi}{25} 4^2 u) \\
	B_{4,1} & 11 & \exp(\frac{i2\pi}{25}(5\cdot 4\cdot1+4^2 u)) \\
	B_{4,2} & 11 & \exp(\frac{i2\pi}{25}(5\cdot 4\cdot2+4^2 u)) \\
	B_{4,3} & 11 & \exp(\frac{i2\pi}{25}(5\cdot 4\cdot3+4^2 u)) \\
	B_{4,4} & 11 & \exp(\frac{i2\pi}{25}(5\cdot 4\cdot4+4^2 u)) \\ 
	\end{array}}
\end{equation}

\section{Modular permutations}

\subsection{Permutatations from $u=1$ to $u=4$.}
$$\scalemath{0.93}{
\begin{array}{c||c|c|c|c|c|c|c|c}
 I_{1,0} & I_{1,0} & I_{1,0} & I_{1,0} & I_{1,0} & I_{1,0} & I_{1,0} & I_{1,0} & I_{1,0} \\
\hline
 I_{1,1} & I_{1,2} & I_{1,2} & I_{1,2} & I_{1,2} & I_{1,3} & I_{1,3} & I_{1,3} & I_{1,3} \\
 I_{1,2} & I_{1,4} & I_{1,4} & I_{1,4} & I_{1,4} & I_{1,1} & I_{1,1} & I_{1,1} & I_{1,1} \\
 I_{1,3} & I_{1,1} & I_{1,1} & I_{1,1} & I_{1,1} & I_{1,4} & I_{1,4} & I_{1,4} & I_{1,4} \\
 I_{1,4} & I_{1,3} & I_{1,3} & I_{1,3} & I_{1,3} & I_{1,2} & I_{1,2} & I_{1,2} & I_{1,2} \\
\hline
 I_{1,5} & I_{1,5} & I_{1,6} & A_{1,0} & A_{2,0} & I_{1,5} & I_{1,6} & A_{1,0} & A_{2,0} \\
 I_{1,6} & I_{1,6} & I_{1,5} & A_{2,0} & A_{1,0} & I_{1,6} & I_{1,5} & A_{2,0} & A_{1,0} \\
\hline
 A_{1,0} & A_{1,0} & A_{2,0} & I_{1,5} & I_{1,6} & A_{1,0} & A_{2,0} & I_{1,5} & I_{1,6} \\
 A_{1,1} & A_{1,1} & A_{2,6} & A_{1,1} & A_{2,6} & A_{1,1} & A_{2,6} & A_{1,1} & A_{2,6} \\
 A_{1,2} & A_{1,2} & A_{2,1} & A_{2,1} & A_{1,2} & A_{1,2} & A_{2,1} & A_{2,1} & A_{1,2} \\
 A_{1,3} & A_{1,3} & A_{2,7} & A_{1,3} & A_{2,7} & A_{1,3} & A_{2,7} & A_{1,3} & A_{2,7} \\
 A_{1,4} & A_{1,4} & A_{2,2} & A_{1,4} & A_{2,2} & A_{1,4} & A_{2,2} & A_{1,4} & A_{2,2} \\
 A_{1,5} & A_{1,5} & A_{2,8} & A_{1,5} & A_{2,8} & A_{1,5} & A_{2,8} & A_{1,5} & A_{2,8} \\
 A_{1,6} & A_{1,6} & A_{2,3} & A_{2,3} & A_{1,6} & A_{1,6} & A_{2,3} & A_{2,3} & A_{1,6} \\
 A_{1,7} & A_{1,7} & A_{2,9} & A_{2,9} & A_{1,7} & A_{1,7} & A_{2,9} & A_{2,9} & A_{1,7} \\
 A_{1,8} & A_{1,8} & A_{2,4} & A_{2,4} & A_{1,8} & A_{1,8} & A_{2,4} & A_{2,4} & A_{1,8} \\
 A_{1,9} & A_{1,9} & A_{2,10} & A_{1,9} & A_{2,10} & A_{1,9} & A_{2,10} & A_{1,9} & A_{2,10} \\
 A_{1,10} & A_{1,10} & A_{2,5} & A_{2,5} & A_{1,10} & A_{1,10} & A_{2,5} & A_{2,5} & A_{1,10} \\
\hline
 A_{2,0} & A_{2,0} & A_{1,0} & I_{1,6} & I_{1,5} & A_{2,0} & A_{1,0} & I_{1,6} & I_{1,5} \\
 A_{2,1} & A_{2,1} & A_{1,2} & A_{1,2} & A_{2,1} & A_{2,1} & A_{1,2} & A_{1,2} & A_{2,1} \\
 A_{2,2} & A_{2,2} & A_{1,4} & A_{2,2} & A_{1,4} & A_{2,2} & A_{1,4} & A_{2,2} & A_{1,4} \\
 A_{2,3} & A_{2,3} & A_{1,6} & A_{1,6} & A_{2,3} & A_{2,3} & A_{1,6} & A_{1,6} & A_{2,3} \\
 A_{2,4} & A_{2,4} & A_{1,8} & A_{1,8} & A_{2,4} & A_{2,4} & A_{1,8} & A_{1,8} & A_{2,4} \\
 A_{2,5} & A_{2,5} & A_{1,10} & A_{1,10} & A_{2,5} & A_{2,5} & A_{1,10} & A_{1,10} & A_{2,5} \\
 A_{2,6} & A_{2,6} & A_{1,1} & A_{2,6} & A_{1,1} & A_{2,6} & A_{1,1} & A_{2,6} & A_{1,1} \\
 A_{2,7} & A_{2,7} & A_{1,3} & A_{2,7} & A_{1,3} & A_{2,7} & A_{1,3} & A_{2,7} & A_{1,3} \\
 A_{2,8} & A_{2,8} & A_{1,5} & A_{2,8} & A_{1,5} & A_{2,8} & A_{1,5} & A_{2,8} & A_{1,5} \\
 A_{2,9} & A_{2,9} & A_{1,7} & A_{1,7} & A_{2,9} & A_{2,9} & A_{1,7} & A_{1,7} & A_{2,9} \\
 A_{2,10} & A_{2,10} & A_{1,9} & A_{2,10} & A_{1,9} & A_{2,10} & A_{1,9} & A_{2,10} & A_{1,9} \\
\hline
 B_{1,0} & B_{3,1} & B_{3,1} & B_{3,1} & B_{3,1} & B_{2,1} & B_{2,1} & B_{2,1} & B_{2,1} \\
 B_{1,1} & B_{3,3} & B_{3,3} & B_{3,3} & B_{3,3} & B_{2,4} & B_{2,4} & B_{2,4} & B_{2,4} \\
 B_{1,2} & B_{3,0} & B_{3,0} & B_{3,0} & B_{3,0} & B_{2,2} & B_{2,2} & B_{2,2} & B_{2,2} \\
 B_{1,3} & B_{3,2} & B_{3,2} & B_{3,2} & B_{3,2} & B_{2,0} & B_{2,0} & B_{2,0} & B_{2,0} \\
 B_{1,4} & B_{3,4} & B_{3,4} & B_{3,4} & B_{3,4} & B_{2,3} & B_{2,3} & B_{2,3} & B_{2,3} \\
\hline
 B_{2,0} & B_{1,0} & B_{1,0} & B_{1,0} & B_{1,0} & B_{4,2} & B_{4,2} & B_{4,2} & B_{4,2} \\
 B_{2,1} & B_{1,2} & B_{1,2} & B_{1,2} & B_{1,2} & B_{4,0} & B_{4,0} & B_{4,0} & B_{4,0} \\
 B_{2,2} & B_{1,4} & B_{1,4} & B_{1,4} & B_{1,4} & B_{4,3} & B_{4,3} & B_{4,3} & B_{4,3} \\
 B_{2,3} & B_{1,1} & B_{1,1} & B_{1,1} & B_{1,1} & B_{4,1} & B_{4,1} & B_{4,1} & B_{4,1} \\
 B_{2,4} & B_{1,3} & B_{1,3} & B_{1,3} & B_{1,3} & B_{4,4} & B_{4,4} & B_{4,4} & B_{4,4} \\
\hline
 B_{3,0} & B_{4,1} & B_{4,1} & B_{4,1} & B_{4,1} & B_{1,1} & B_{1,1} & B_{1,1} & B_{1,1} \\
 B_{3,1} & B_{4,3} & B_{4,3} & B_{4,3} & B_{4,3} & B_{1,4} & B_{1,4} & B_{1,4} & B_{1,4} \\
 B_{3,2} & B_{4,0} & B_{4,0} & B_{4,0} & B_{4,0} & B_{1,2} & B_{1,2} & B_{1,2} & B_{1,2} \\
 B_{3,3} & B_{4,2} & B_{4,2} & B_{4,2} & B_{4,2} & B_{1,0} & B_{1,0} & B_{1,0} & B_{1,0} \\
 B_{3,4} & B_{4,4} & B_{4,4} & B_{4,4} & B_{4,4} & B_{1,3} & B_{1,3} & B_{1,3} & B_{1,3} \\
\hline
 B_{4,0} & B_{2,0} & B_{2,0} & B_{2,0} & B_{2,0} & B_{3,2} & B_{3,2} & B_{3,2} & B_{3,2} \\
 B_{4,1} & B_{2,2} & B_{2,2} & B_{2,2} & B_{2,2} & B_{3,0} & B_{3,0} & B_{3,0} & B_{3,0} \\
 B_{4,2} & B_{2,4} & B_{2,4} & B_{2,4} & B_{2,4} & B_{3,3} & B_{3,3} & B_{3,3} & B_{3,3} \\
 B_{4,3} & B_{2,1} & B_{2,1} & B_{2,1} & B_{2,1} & B_{3,1} & B_{3,1} & B_{3,1} & B_{3,1} \\
 B_{4,4} & B_{2,3} & B_{2,3} & B_{2,3} & B_{2,3} & B_{3,4} & B_{3,4} & B_{3,4} & B_{3,4} 
\end{array}}
$$

\subsection{Permutations from $u=2$ to $u=3$}
\par
$$\scalemath{0.93}{\begin{array}{c||c|c|c|c|c|c|c|c}
 I_{1,0} & I_{1,0} & I_{1,0} & I_{1,0} & I_{1,0} & I_{1,0} & I_{1,0} & I_{1,0} & I_{1,0} \\
\hline
 I_{1,1} & I_{1,2} & I_{1,2} & I_{1,2} & I_{1,2} & I_{1,3} & I_{1,3} & I_{1,3} & I_{1,3} \\
 I_{1,2} & I_{1,4} & I_{1,4} & I_{1,4} & I_{1,4} & I_{1,1} & I_{1,1} & I_{1,1} & I_{1,1} \\
 I_{1,3} & I_{1,1} & I_{1,1} & I_{1,1} & I_{1,1} & I_{1,4} & I_{1,4} & I_{1,4} & I_{1,4} \\
 I_{1,4} & I_{1,3} & I_{1,3} & I_{1,3} & I_{1,3} & I_{1,2} & I_{1,2} & I_{1,2} & I_{1,2} \\
\hline
 I_{1,5} & I_{1,5} & I_{1,6} & A_{1,0} & A_{2,0} & I_{1,5} & I_{1,6} & A_{1,0} & A_{2,0} \\
 I_{1,6} & I_{1,6} & I_{1,5} & A_{2,0} & A_{1,0} & I_{1,6} & I_{1,5} & A_{2,0} & A_{1,0} \\
\hline
 A_{1,0} & A_{1,0} & A_{2,0} & I_{1,5} & I_{1,6} & A_{1,0} & A_{2,0} & I_{1,5} & I_{1,6} \\
 A_{1,1} & A_{1,1} & A_{2,6} & A_{1,1} & A_{2,6} & A_{1,1} & A_{2,6} & A_{1,1} & A_{2,6} \\
 A_{1,2} & A_{1,2} & A_{2,1} & A_{2,1} & A_{1,2} & A_{1,2} & A_{2,1} & A_{2,1} & A_{1,2} \\
 A_{1,3} & A_{1,3} & A_{2,7} & A_{1,3} & A_{2,7} & A_{1,3} & A_{2,7} & A_{1,3} & A_{2,7} \\
 A_{1,4} & A_{1,4} & A_{2,2} & A_{1,4} & A_{2,2} & A_{1,4} & A_{2,2} & A_{1,4} & A_{2,2} \\
 A_{1,5} & A_{1,5} & A_{2,8} & A_{1,5} & A_{2,8} & A_{1,5} & A_{2,8} & A_{1,5} & A_{2,8} \\
 A_{1,6} & A_{1,6} & A_{2,3} & A_{2,3} & A_{1,6} & A_{1,6} & A_{2,3} & A_{2,3} & A_{1,6} \\
 A_{1,7} & A_{1,7} & A_{2,9} & A_{2,9} & A_{1,7} & A_{1,7} & A_{2,9} & A_{2,9} & A_{1,7} \\
 A_{1,8} & A_{1,8} & A_{2,4} & A_{2,4} & A_{1,8} & A_{1,8} & A_{2,4} & A_{2,4} & A_{1,8} \\
 A_{1,9} & A_{1,9} & A_{2,10} & A_{1,9} & A_{2,10} & A_{1,9} & A_{2,10} & A_{1,9} & A_{2,10} \\
 A_{1,10} & A_{1,10} & A_{2,5} & A_{2,5} & A_{1,10} & A_{1,10} & A_{2,5} & A_{2,5} & A_{1,10} \\
\hline
 A_{2,0} & A_{2,0} & A_{1,0} & I_{1,6} & I_{1,5} & A_{2,0} & A_{1,0} & I_{1,6} & I_{1,5} \\
 A_{2,1} & A_{2,1} & A_{1,2} & A_{1,2} & A_{2,1} & A_{2,1} & A_{1,2} & A_{1,2} & A_{2,1} \\
 A_{2,2} & A_{2,2} & A_{1,4} & A_{2,2} & A_{1,4} & A_{2,2} & A_{1,4} & A_{2,2} & A_{1,4} \\
 A_{2,3} & A_{2,3} & A_{1,6} & A_{1,6} & A_{2,3} & A_{2,3} & A_{1,6} & A_{1,6} & A_{2,3} \\
 A_{2,4} & A_{2,4} & A_{1,8} & A_{1,8} & A_{2,4} & A_{2,4} & A_{1,8} & A_{1,8} & A_{2,4} \\
 A_{2,5} & A_{2,5} & A_{1,10} & A_{1,10} & A_{2,5} & A_{2,5} & A_{1,10} & A_{1,10} & A_{2,5} \\
 A_{2,6} & A_{2,6} & A_{1,1} & A_{2,6} & A_{1,1} & A_{2,6} & A_{1,1} & A_{2,6} & A_{1,1} \\
 A_{2,7} & A_{2,7} & A_{1,3} & A_{2,7} & A_{1,3} & A_{2,7} & A_{1,3} & A_{2,7} & A_{1,3} \\
 A_{2,8} & A_{2,8} & A_{1,5} & A_{2,8} & A_{1,5} & A_{2,8} & A_{1,5} & A_{2,8} & A_{1,5} \\
 A_{2,9} & A_{2,9} & A_{1,7} & A_{1,7} & A_{2,9} & A_{2,9} & A_{1,7} & A_{1,7} & A_{2,9} \\
 A_{2,10} & A_{2,10} & A_{1,9} & A_{2,10} & A_{1,9} & A_{2,10} & A_{1,9} & A_{2,10} & A_{1,9} \\
\hline
 B_{1,0} & B_{3,0} & B_{3,0} & B_{3,0} & B_{3,0} & B_{2,4} & B_{2,4} & B_{2,4} & B_{2,4} \\
 B_{1,1} & B_{3,2} & B_{3,2} & B_{3,2} & B_{3,2} & B_{2,2} & B_{2,2} & B_{2,2} & B_{2,2} \\
 B_{1,2} & B_{3,4} & B_{3,4} & B_{3,4} & B_{3,4} & B_{2,0} & B_{2,0} & B_{2,0} & B_{2,0} \\
 B_{1,3} & B_{3,1} & B_{3,1} & B_{3,1} & B_{3,1} & B_{2,3} & B_{2,3} & B_{2,3} & B_{2,3} \\
 B_{1,4} & B_{3,3} & B_{3,3} & B_{3,3} & B_{3,3} & B_{2,1} & B_{2,1} & B_{2,1} & B_{2,1} \\
\hline
 B_{2,0} & B_{1,1} & B_{1,1} & B_{1,1} & B_{1,1} & B_{4,3} & B_{4,3} & B_{4,3} & B_{4,3} \\
 B_{2,1} & B_{1,3} & B_{1,3} & B_{1,3} & B_{1,3} & B_{4,1} & B_{4,1} & B_{4,1} & B_{4,1} \\
 B_{2,2} & B_{1,0} & B_{1,0} & B_{1,0} & B_{1,0} & B_{4,4} & B_{4,4} & B_{4,4} & B_{4,4} \\
 B_{2,3} & B_{1,2} & B_{1,2} & B_{1,2} & B_{1,2} & B_{4,2} & B_{4,2} & B_{4,2} & B_{4,2} \\
 B_{2,4} & B_{1,4} & B_{1,4} & B_{1,4} & B_{1,4} & B_{4,0} & B_{4,0} & B_{4,0} & B_{4,0} \\
\hline
 B_{3,0} & B_{4,1} & B_{4,1} & B_{4,1} & B_{4,1} & B_{1,3} & B_{1,3} & B_{1,3} & B_{1,3} \\
 B_{3,1} & B_{4,3} & B_{4,3} & B_{4,3} & B_{4,3} & B_{1,1} & B_{1,1} & B_{1,1} & B_{1,1} \\
 B_{3,2} & B_{4,0} & B_{4,0} & B_{4,0} & B_{4,0} & B_{1,4} & B_{1,4} & B_{1,4} & B_{1,4} \\
 B_{3,3} & B_{4,2} & B_{4,2} & B_{4,2} & B_{4,2} & B_{1,2} & B_{1,2} & B_{1,2} & B_{1,2} \\
 B_{3,4} & B_{4,4} & B_{4,4} & B_{4,4} & B_{4,4} & B_{1,0} & B_{1,0} & B_{1,0} & B_{1,0} \\
\hline
 B_{4,0} & B_{2,2} & B_{2,2} & B_{2,2} & B_{2,2} & B_{3,2} & B_{3,2} & B_{3,2} & B_{3,2} \\
 B_{4,1} & B_{2,4} & B_{2,4} & B_{2,4} & B_{2,4} & B_{3,0} & B_{3,0} & B_{3,0} & B_{3,0} \\
 B_{4,2} & B_{2,1} & B_{2,1} & B_{2,1} & B_{2,1} & B_{3,3} & B_{3,3} & B_{3,3} & B_{3,3} \\
 B_{4,3} & B_{2,3} & B_{2,3} & B_{2,3} & B_{2,3} & B_{3,1} & B_{3,1} & B_{3,1} & B_{3,1} \\
 B_{4,4} & B_{2,0} & B_{2,0} & B_{2,0} & B_{2,0} & B_{3,4} & B_{3,4} & B_{3,4} & B_{3,4} 
\end{array}}$$

\end{document}